\documentclass[11pt]{amsart}
\usepackage{amssymb} 
\usepackage{color} 
\newtheorem{theorem}{Theorem} [section]
\newtheorem{lemma}{Lemma}[section]

\newtheorem{conjecture}{Conjecture}[section]
\newtheorem{remark}{Remark}[section]
\newtheorem{question}{Question}[section]
\begin{document}
\title{On existence of PI-exponent of algebras with involution}
\author[D.D Repov\v s and M.V. Zaicev]
{Du\v san D. Repov\v s and Mikhail V. Zaicev}
\address{Du\v san D. Repov\v s \\Faculty of Education, and
Faculty of Mathematics and Physics, University of Ljubljana
\& Institute of Mathematics and Physics, Ljubljana, 1000, Slovenia}
\email{dusan.repovs@guest.arnes.si}
\address{Mikhail V. Zaicev \\Department of Algebra\\ Faculty of Mathematics and
Mechanics\\  Moscow State University \\ Moscow,119992, Russia\\
Moscow Center of Fundamental and Applied Mathematics, Moscow, 119991 Russia}
\email{zaicevmv@mail.ru}
\begin{abstract}
We study polynomial identities of algebras with involution of nonassociative algebras over a field
of characteristic zero. We prove that the growth of the sequence of $*$-codimensions of a finite-dimensional 
algebra is exponentially bounded. We construct a series of finite-dimensional algebras
with fractional $*$-PI-exponent. We also construct a family of infinite-dimensional algebras
$C_\alpha$ such that ${\rm exp}^*(C_\alpha)$ does not exist.
\end{abstract}   
\date{}
\maketitle

\section{Introduction}\label{s1}
Let $A$ be an algebra over a field $\Phi$ of characteristic zero. One of the modern approaches to the study 
of polynomial identities of $A$ is to investigate  their numerical invariants. The most important
numerical characteristic of identities of $A$ is the sequence $\{c_n(A)\}$ of codimensions and its 
asymptotic behavior. For a wide class of algebras, the growth of the sequence $\{c_n(A)\}$ is exponentially bounded.
This class includes associative PI-algebras \cite{R,L},  finite-dimensional algebras
of arbitrary signature \cite{BD,GZ}, affine Kac-Moody algebras \cite{Z}, infinite-dimensional
simple Lie algebras of Cartan type \cite{M}, Virasoro algebra, Novikov algebras \cite{Dz}, and many
others.
 
In the case of exponential upper bound, the corresponding  sequence of roots $\{\sqrt[n]{c_n(A)}\}$ is
bounded and its lower and upper limits
$$
\underline{{\rm exp}}(A)=\liminf_{n\to\infty}\sqrt[n]{c_n(A)},\quad
\overline{{\rm exp}}(A)=\limsup_{n\to\infty}\sqrt[n]{c_n(A)}
$$
are called the {\it lower} and the {\it upper} PI-{\it exponent} of $A$, respectively. In the case when $\underline{{\rm exp}}(A)=
\overline{{\rm exp}}(A)$,  the ordinary limit
$$
{\rm exp}(A)=\lim_{n\to\infty}\sqrt[n]{c_n(A)}
$$
is called the ({\it ordinary}) PI-{\it exponent} of $A$.

In the late 1980's, S. Amitsur conjectured that the PI-exponent of any associative PI-algebra 
exists and is a nonnegative integer. Amitsur's conjecture was confirmed in \cite{GZ1}. It was also
proved for finite-dimensional Lie algebras \cite{Z1},  Jordan algebras \cite{GSZ}, and some
others. The class of algebras for which Amitsur's conjecture was partially confirmed is much wider. Namely,
the existence (but not the  integrality, in general) was proved in a series of papers.

 For example, it 
was shown in \cite{GZ2} that the PI-exponent exists for any finite-dimensional simple algebra.
The question about existence of PI-expo\-nents is one of the main problems of numerical theory of
polynomial identities. Until now, only two results about algebras without PI-exponent have been proved.
An example of a
two-step left-nilpotent algebra without PI-exponent  was
constructed in \cite{ERA}.
Analogous  result for unitary algebras was obtained in \cite{RZ}.

If an algebra $A$ is equipped with an additional structure (like an involution or a group grading),
then one may consider identities with involution, graded identities, etc. Recall that in the
associative case, the celebrated theorem of Amitsur \cite{A} states that if $A$ is an algebra with
involution $*: A\to A,$ satisfying a $*$-polynomial identity, then $A$ satisfies an ordinary (non-involution)
polynomial identity. As a consequence, the sequence of $*$-codimensions $\{c_n^*(A)\}$ is exponentially
bounded. In  \cite{GZ3, GPV} the existence and integrality of ${\rm exp}^*(A)$ was proved
for any associative PI-algebra with involution.

In the present paper we shall show that the class of algebras with exponentially bounded $*$-codimension
sequence is sufficiently large. In particular, it contains all finite-dimensional algebras.

{\bf Theorem A} (see Theorem~\ref{t1} in Section~\ref{s3}).
%\begin{theorem}\label{t1}
{\it Let $A$ be a finite-dimensional algebra with involution $*\colon A\to A$ and $d=\dim A$. Then $*$-codimensions
of $A$ satisfy the following inequality
$$
c_n^*(A)\le d^{n+1}.
$$
}

 Nevertheless,
as it will be shown, the results of \cite{GZ3, GPV} cannot be generalized 
to the general nonassociative case. We shall construct a series of finite-dimensional algebras with
fractional $*$-PI-exponent.   For any integer $T\ge 2$ we shall construct an algebra $A_T$ with the following
property.

{\bf Theorem B} (see Theorem~\ref{t2} in Section~\ref{s4}).
{\it The $*$-PI-exponent of algebra $A_T$ exists and 
$$
exp^*(A_T) = \frac{1}{\theta_T^{\theta_T}(1-\theta_T)^{1-\theta_T}},
$$
where $\theta_T=\frac{1}{2T+1}$.
}

We shall also present a family of algebras $C_\alpha$ with involution $*$
which has an exponentially bounded sequence $\{c_n^*(C_\alpha)\}$ such that ${\rm exp}^*(C_\alpha)$ does
not exist.

{\bf Theorem C} (see Theorem~\ref{t3} in Section~\ref{s5}).
{\it For any real number $\alpha>1$ there exists an algebra $C_\alpha$ such that
$$
\underline{exp}^*(C_\alpha)=1,~~\overline{exp}^*(C_\alpha)=\alpha.
$$
}

The necessary background on numerical theory of polynomial identities can be found in \cite{GZBook}.

\section{Preliminaries}\label{s2}

Let $A$ be an algebra with involution $*\colon A\to A$ over a field $\Phi$ of $\rm{char}~\Phi=0$. Recall
that an element $a\in A$ is called {\it symmetric} if $a^*=a$, whereas an element $b\in A$ is called {\it
skew-symmetric} if $b^*=-b$. Denote
$$
A^+=\{a\in A|a^*=a\},~~ A^-=\{b\in A|b^*=-b\}.
$$
Obviously, we have a vector space decomposition $A=A^+\oplus A^-$. In order to study $*$-polynomial identities
 we need to introduce free objects in the following way.
 
Let $\Phi\{X,Y\}$ be a free (nonassociative) algebra over $\Phi$ with the set of free generators
$X\cup Y$, $X=\{x_1,x_2,\ldots\}, Y=\{y_1,y_2,\ldots\}$. A map $*:X\cup Y\to X\cup Y$ such that
$x^*_i=x_i,y^*_i=-y_i, i=1,2,\ldots$, can be naturaly extended to an involution on $\Phi\{X,Y\}$.
A polynomial $f=f(x_1,\dots,x_m,y_1,\ldots,y_n)\in\Phi\{X,Y\}$ is said to be {\it a $*$-identity}
 of $A$ if $$f(a_1,\ldots,a_m,b_1,\ldots,b_n)=0,
 \
 \hbox{for all}
 \
 a_1,\ldots a_m\in A^+,
b_1,\ldots,b_n\in A^-.$$ 

Denote by $Id^*(A)$ the set of all $*$-identities of $A$ in $\Phi\{X,Y\}$.
Then  $Id^*(A)$ is an ideal of $\Phi\{X,Y\}$ and it is  stable under involution $*$ and endomorphisms compatible
with $*$.

Given $0\le k\le n$, denote the space of all multilinear polynomials in $\Phi\{X,Y\}$
in $k$ symmetric variables $x_1,\ldots,x_k$ and  $n-k$ skew-symmetric variables $y_1,\ldots,y_{n-k}$  by $P^*_{k,n-k}$. Denote
also $$P_n^*=P_{0,n}^*\oplus P_{1,n-1}^*\oplus\cdots\oplus P_{n,0}^*.$$ Clearly, the intersection
$P_{k,n-k}^*\cap Id^*(A)$ is the subspace of all multilinear  $*$-identities of $A$ in $k$ symmetric
and $n-k$ skew-symmetric variables.

The following value
$$
c^*_{k,n-k}(A)=\dim\frac{P^*_{k,n-k}}{P^*_{k,n-k}\cap Id^*(A)}
$$
is called the {\it partial} $(k,n-k)~ $ $*$-{\it codimension} of $A$, whereas the value
$$
c^*_n(A)=\sum_{k=0}^n {n\choose k} c^*_{k,n-k}(A)
$$  
is called the ({\it total})~~ $*$-{\it codimension} of $A$. We shall also use the following notations
$$
P^*_{k,n-k}(A)=\frac{P^*_{k,n-k}}{P^*_{k,n-k}\cap Id^*(A)},~~
P^*_{n}(A)=\frac{P^*_{n}}{P^*_{n}\cap Id^*(A)}.
$$

\section{$*$-codimensions of finite-dimensional algebras}\label{s3}

Let $A$ be a finite-dimensional algebra with involution $*\colon A\to A$, where $\dim A=d$. Recall that
$A^+$ and $A^-$ are the subspaces of symmetric and skew-symmetric elements of $A$, respectively. In order to get an
exponential upper bound for $c^*_n(A),$ we shall follow the approach of \cite{BD}. Choose a basis $a_1,\ldots,a_p$
of $A^+$ and a basis $b_1,\ldots,b_q$ of $A^-$. If $f(x_1,\ldots,x_k,y_1,\ldots,y_{n-k})\in P_{k,n-k}^*$
is a multilinear $*$-polynomial in $k$ symmetric variables $x_1,\ldots,x_k$ and  $n-k$ skew-symmetric variables
$y_1,\ldots,y_{n-k},$ then $f$ is a $*$-identity of $A$ if and only if $\varphi(f)=0,$ for all evaluations
$\varphi$ such that
\begin{equation}\label{ef1}
\varphi(x_i)\in \{a_1,\ldots,a_p\},~1\le i\le k,~~ \varphi(y_j)\in \{b_1,\ldots,b_q\},~1\le j\le n-k.
\end{equation}

Denote $N=\dim P_{k,n-k}^*$. Fix a basis $g_1,\ldots,g_N$ of $P_{k,n-k}^*$ and write $f$ as a linear combination
$f=\alpha_1g_1+\cdots+\alpha_Ng_N$. Then the value $\varphi(f)$ for $\varphi$ of the type (\ref{ef1}) can be
written as
$$
\varphi(f)=\lambda_1a_1+\cdots+\lambda_pa_p+\mu_1b_1+\cdots+b_q\mu_q,
$$
where all $\lambda_1,\ldots,\lambda_p,\mu_1,\ldots,\mu_q$ are linear combinations of $\alpha_1,\ldots,\alpha_N$.
Hence $\varphi(f)=0$ if and only if
\begin{equation}\label{ef2}
\lambda_1=\cdots=\lambda_p=\mu_1=\cdots\mu_q=0.
\end{equation}
The total number of evaluations $\varphi$ of  type (\ref{ef1}) is equal to $p^kq^{n-k}$. It follows that $f\equiv 0$
is a $*$-identity of $A$ if and only if the $N$-tuple $(\alpha_1,\ldots,\alpha_N)$ is the solution of system $S$ of
$p^kq^{n-k}(p+q)$ linear equations of type (\ref{ef2}). 

Denote by $U$ the subspace of all solutions of system $S$
in the space $V$ of all $N$-tuples $(\alpha_1,\ldots,\alpha_N)$. Then $\dim U=N-r$, where $r={\rm rank}~S$ is the rank of $S$.
Clearly,
\begin{equation}\label{ef3}
r\le p^kq^{n-k}(p+q).
\end{equation}

Since
$$
c_{k,n-k}^*(A)=\rm{codim}_V(U)=r,
$$
it follows from (\ref{ef3}) that $$c_{k,n-k}^*(A)\le p^kq^{n-k}(p+q)$$ and
$$
c^*_n(A)=\sum_{k=0}^n {n\choose k} c_{k,n-k}^*(A) \le (p+q)\sum_{k=0}^n {n\choose k}p^kq^{n-k}=(p+q)^{n+1}.
$$

Recall that $p+q=d=\dim A$. Hence we have proved the first main result of this paper.

\begin{theorem}\label{t1}
Let $A$ be a finite-dimensional algebra with involution $*\colon A\to A$ and $d=\dim A$. Then $*$-codimensions
of $A$ satisfy the following inequality
$$
c_n^*(A)\le d^{n+1}. 
$$
\end{theorem} 
\hfill $\Box$

In the case of exponentially bounded sequence $\{c_n^*(A)\}$, the following natural question arises.
\begin{question}
Does the $*$-PI-exponent
$$
{\rm exp}^*(A)=\lim_{n\to\infty}\sqrt[n]{c_n^*(A)}
$$
exist and what are its possible values?
\end{question}

 In Section~\ref{s1} we mentioned that $c^*_n(A)$ exists and is
a nonnegative integer for any associative $*$-PI-algebra $A$. The following hypotheses look very
natural.
\vskip.1in

\begin{conjecture}\label{c1}
For any finite-dimensional algebra $A$ with involution $*$, its $*$-PI-exponent ${\rm exp}^*(A)$  exists.
 \end{conjecture}

In the light of results of \cite{GMZ}, we can assume that $*$-PI-exponent may take on all real values $\ge 1.$ 

\begin{conjecture}\label{c2}
For any real value $\alpha\ge 1$, there exists an algebra $A_\alpha$ with involution such that 
$*$-PI-exponent of $A_\alpha$ exists and ${\rm exp}^*(A_\alpha)=\alpha$.
\end{conjecture}

\section{Algebras with fractional $*$-PI-exponent}\label{s4}

In this section  we shall discuss $*$-codimension growth of algebras $A_T$ introduced in \cite{SZ}.
We shall prove the existence of $*$-PI-exponents of $A_T$ and compute the precise value   of ${\rm exp}^*(A_T)$. In
Section~\ref{s5} we shall use the properties of $A_T$ for constructing several counterexamples.
 
Recall the structure of $A_T$. Given an integer $T\ge 2$, denote by $A_T$ the algebra with basis
$\{a,b,z_1,\ldots,z_{2T+1}\}$ and with multiplication
$$
z_ia=az_i=z_{i+1}, 1\le i\le 2t,~z_{2T+1}b=bz_{2T+1}=z_1,
$$
where all remaining products are zero. Involution $*:A_T\to A_T$ is defined by
$$
a^*=-a,b^*=b,z_i^*=(-1)^{i+1}z_i
$$
and then $$A^+=<b,z_1,z_3,\ldots,z_{2T+1}>, A^-=<a,z_2,z_4,\ldots,z_{2T}>.$$

We shall need the following two results from \cite{SZ}.

\begin{lemma}\label{l1} (\cite[Lemma 3.7]{SZ})
The 
$*$-codimensions of $A_T$ satisfy the inequality $c^*_n(A_T) \le n^3$, provided that $n\le 2T$.
\end{lemma}

\begin{lemma}\label{l2} (\cite[Corollary 3.8]{SZ})
Let $f\equiv 0$ be a multilinear $*$-identity of $A_T$ of degree $n\le 2T$. Then $f$ is also
an identity of $A_{T+1}$.
\end{lemma}

Note that algebras $A_T$ are commutative and {\it metabelian}, i.e. they satisfy the following identity
$$
(xy)(zt)\equiv 0.
$$
Hence any product of elements $c_1,\ldots,c_n\in A$ can be written in the left-normed form. We shall
omit brackets in the left-normed products, i.e. we shall write $c_1c_2\cdots c_n$ instead of
$(\ldots(c_1c_2)\ldots)c_n$.
 
First, we shall find a lower bound for $*$-codimensions.

\begin{lemma}\label {l3}
The following inequality holds
for all
$
n\ge 2T+2,$
\begin{equation}\label{efr0}
c^*_n(A_T)\ge\frac{1}{n^2} \left(\frac{1}{\theta_T^{\theta_T}(1-\theta_T)^{1-\theta_T}}\right)^{n-2T-1},
\end{equation}
 where 
$$
\theta_T=\frac{1}{2T+1}.
$$
\end{lemma}

{\em Proof.} Write $n$ in the form $n=(2T+1)k+t+1$, where $0\le t \le 2T$. Then the following product
of $n$ basis elements is nonzero
$$
z_1\underbrace{a^{2T}b\cdots a^{2T}b}_ka^t=z_{t+1}\ne 0.
$$ 
Here, we use the  notation $xa^m$ for $x\underbrace{a\cdots a}_m$. Hence  the polynomial
$$
x_0y_1\cdots y_{2T}x_1\cdots y_{2t(k-1)+1}\cdots y_{2Tk}x_k y_{2Tk+1}\cdots y_{2Tk+t}
$$
is not an identity of $A_T$, that is, 
$$P^*_{k+1,2Tk+t}(A_T)\ne 0, ~~
c^*_{k+1,2Tk+t}\ge 1.$$

In particular,
\begin{equation}\label{efr1}
c_n^*(A_T)\ge{n\choose k+1}\ge{n_0\choose k+1}\ge{n_0\choose k},
\end{equation}
where $n=2Tk+k+t+1,n_0=2Tk+k$. 

Using the Stirling formula for factorials we get
\begin{equation}\label{efr1a}
{(2T+1)k\choose k}>\frac{1}{n^2}\frac{((2T+1)k)^{(2T+1)k}}{k^k (2Tk)^{2Tk}}
\end{equation}
$$
=\frac{1}{n^2}
\left(
\frac{1}{\left(\frac{1}{2T+1}\right)^\frac{1}{2T+1}\left(\frac{2T}{2T+1}\right)^\frac{2T}{2T+1}}
\right)^{(2T+1)k}
=\frac{1}{n^2}\left( \frac{1}{\theta_T^{\theta_T}(1-\theta_T)^{1-\theta_T}}\right)^{n_0}
$$
$$
\ge \frac{1}{n^2}\left( \frac{1}{\theta_T^{\theta_T}(1-\theta_T)^{1-\theta_T}}\right)^{n-2T-1},
$$
where $\theta_T=\frac{1}{2T+1}$. 

Finally, combining (\ref{efr1}) and (\ref{efr1a}), we obtain the desired
inequality (\ref{efr0}).
\hfill $\Box$
 
Next, we shall find an upper bound for $c_n^*(A_T)$. First, we restrict the number of nonzero
components $P^*_{k,n-k}(A_T)$ for a fixed  $n$.

\begin{lemma}\label{l4}
Given a positive integer $n$, there are at most three integers $k$, $0\le k\le n,$ such that $P^*_{k,n-k}(A_T)\ne 0$.
Moreover, if
$P^*_{k,n-k}(A_T)\ne 0,$ then
 $$\frac{k-2}{n}\le\frac{1}{2T+1}.$$
\end{lemma}

{\em Proof}.  Clearly, all nonzero products of the basis elements of $A_T$ are of the form
\begin{equation}\label{efr2}
W=z_{2T+1-i}a^ib\underbrace{a^{2T}b\cdots a^{2T}b}_p a^j.
\end{equation}
The number of symmetric factors $k$ is equal to $p+1$ if $i$ is odd, and $k=p+2$ if $i$ is even. The total
number of factors in $W$ is equal to $n=(2T+1)p+i+j+2$. Moreover, $i$ and $j$ in (\ref{efr2})
satisfy inequalities $0\le i,j\le 2T$.  Hence
\begin{equation}\label{efr3}
n-4T-2\le (2T+1)p \le n-2.
\end{equation}

Clearly, there are at most two integers $p$ satisfying (\ref{efr3}). Since $k=p+1$ or $p+2$, at 
most 3 components $P^*_{k,n-k}(A_T)$ can be nonzero. Finally, according to (\ref{efr3}), we
have
$$
\frac{k-2}{n}\le \frac{p}{n}\le\frac{n-2}{(2T+1)n}\le\frac{1}{2T+1}.
$$
\hfill $\Box$

\begin{lemma}\label{l6}
Let $n\le 2T+2$. Then $c^*_{k,n-k}(A_T)\le (2T+1)^3$.
\end{lemma}
{\em Proof}. As it was mentioned earlier, all nonzero products of the basis elements of $A_T$ 
are of the form
$$ 
z_ja^pb\underbrace{a^{2T}b\cdots a^{2T}b}_ka^q,~~1\le j\le 2T+1,~~ 0\le p,q\le 2T.
$$
Hence all nonzero modulo $Id^*(A_T)$ multilinear monomials are of the form
\begin{equation}\label{efr5}
wy_{\sigma(1)}\cdots y_{\sigma(p)}x_{\tau(1)}y_{\sigma(p+1)}\cdots y_{\sigma(p+2T)}x_{\tau(2)}
\cdots
\end{equation} 
$$
y_{\sigma(2Tk-2T+p+1)}\cdots y_{\sigma(2Tk+p)} x_{\tau(k+1)} y_{\sigma(2Tk+p+1)}\cdots
y_{\sigma(2Tk+p+q)},,
$$
where $\sigma\in S_{2Tk+p+q},\tau\in S_{k+1}$, and $w$ is either $x_0$ or $y_0$.

Moreover, any monomial (\ref{efr5}) coincides (modulo $Id^*(A_T)$) with the special case (\ref{efr5}) when
$\sigma=1,\tau=1$. Hence, we have at most $(2T+1)^3$ linearly independent elements in
$P^*_{k,n-k}(A_T)$, and so we are done.
\hfill $\Box$

\begin{lemma}\label{l7}
For all $n\ge 2T+2$, we have
$$
c^*_n(A_T) \le 3(2T+1)^3n^3\left(\frac{1}{\theta_T^{\theta_T}(1-\theta_T)^{1-\theta_T}}\right)^n.
$$
\end{lemma}
{\em Proof.} First we compute an upper bound for $c^*_{k,n-k}(A_T)$, provided that $P^*_{k,n-k}(A_T)\ne 0$.
Note that
$$
{n\choose k}\le n^2{n\choose k-2} \le n^3\frac{n^n}{m^m(n-m)^{n-m}},
$$ 
by the Stirling formula, where $m=k-2$. 

Since the function $$\frac{1}{x^x(1-x)^{1-x}}$$ is nondecreasing on
$(0,\frac{1}{2})$, we have by Lemma \ref{l4},
\begin{equation}\label{efr7}
{n\choose k}\le n^3 \left(\frac{1}{(m/n)^{m/n}(1-m/n)^{1-m/n}} \right)^n
\le n^3\left(\frac{1}{\theta_T^{\theta_T}(1-\theta_T)^{1-\theta_T}}\right)^n.
\end{equation}

Now relation (\ref{efr7}), Lemma \ref{l4}, and Lemma \ref{l6} imply
$$
c^*_n(A_T)=\sum_{k=0}^n {n\choose k}c^*_{k,n-k}(A_T)\le 3(2T+1)^3n^3 
\left(\frac{1}{\theta_T^{\theta_T}(1-\theta_T)^{1-\theta_T}}\right)^n.
$$
\hfill $\Box$

Finally, Lemma \ref{l3} and Lemma \ref{l7} imply the second main result of this paper.

\begin{theorem}\label{t2}
The
$*$-PI-exponent of algebra $A_T$ exists and 
$$
{\rm exp}^*(A_T) = \frac{1}{\theta_T^{\theta_T}(1-\theta_T)^{1-\theta_T}},
$$
where $\theta_T=\frac{1}{2T+1}$.
\end{theorem}
\hfill $\Box$

\section{Algebras without $*$-PI-exponent}\label{s5}

We modify construction of the algebra from Section~\ref{s4}. Denote by  $\widetilde A_T$ an 
infinite-dimensional algebra with the basis
$$
a,b_i, z^i_j,~~ 1\le j\le 2T+1,~~  i=1,2,\ldots
$$ and multiplication table
$$
az^i_j=z^i_ja=z^i_{j+1},1\le j\le 2T, ~~
  b_iz^i_{2T+1}=z^i_{2T+1}b_i=z^{i+1}_1.
$$

Involution $*:\widetilde A_T\to\widetilde A_T$ is defined as follows
$$
a^*=-a, ~~ b_i^*=b_i, ~~ (z^i_j)^*=(-1)^{j+1}z^i_j,~~ 1\le j\le 2T+1,~~ i=1,2,\ldots~~.
$$

\begin{lemma}\label{l8}
A multilinear polynomial $f\in  P^*_{k,n-k}$ of degree $n\le 2T$ is a $*$-identity
of $\widetilde A_T$ if and only if $f$ is a $*$-identity of $A_T$.
\end{lemma}
{\em Proof.} First, note that $P^*_{k,n-k}(A_T)=P^*_{k,n-k}(\widetilde A_T)=0$, when
 $n\le 2T$ and $3\le k\le n$. 

Let $k=0$. Then both $A_T$ and $\widetilde A_T$ satisfy the following identity
$$
y_{t+1}y_{\sigma(1)}\cdots y_{\sigma(t)}=y_{t+1}y_1\cdots y_t,
$$
for any $\sigma\in S_t$ and $t\le 2T-1$. Hence, modulo $Id^*(A_T)$ (and modulo $Id^*(\widetilde A_T)$),
the polynomial $f$ coincides with linear combination $$f=\lambda_2w_2+\cdots+\lambda_nw_n,
\
\hbox{ where}
\
w_j=y_jy_1\cdots y_{j-1}y_{j+1}\cdots y_n.$$

 Let for example, $\lambda_n\ne 0$. Then $\varphi(f)\ne 0$
in $A_T$ and $\widetilde\varphi(f)\ne 0$ in $\widetilde A_T$ for evaluations $\varphi,\widetilde\varphi$, 
where
$$
\varphi(y_n)=z_1,\varphi(y_j)=a~ {\rm in}~A_T, 2\le j\le n-1,~\widetilde\varphi(y_n)=z^1_1, 
\widetilde\varphi(y_j)=a ~{\rm in}~\widetilde A_T, 2\le j\le n-1.
$$

Now let $k=1$. Then all monomials $y_1\cdots y_jx_1y_{j+1}\cdots y_t$ are identities of $A_T$ and
$\widetilde A_T$ if $3\le j\le t\le n-1$. Since
$$
x_1y_{\sigma(1)}\cdots y_{\sigma(n-1)}\equiv x_1y_1\cdots y_{n-1}, ~~
{\rm for~all}~~ \sigma\in S_{n-1}
$$
${\rm mod}~Id^*(A_T)$ and ${\rm mod}~Id^*(\widetilde A_T)$, it follows that  $f=\lambda x_1y_1\cdots y_{n-1}$, with
$0\ne\lambda\in\Phi$. Hence $f\not\in Id^*(A_T)$ and $f\not\in Id^*(\widetilde A_T)$.

Finally, let $k=2$. Then modulo  $Id^*(A_T)$ and modulo $Id^*(\widetilde A_T)$, any multilinear
$*$-polynomial is a linear combination of monomials
$$
w_p=x_1y_1\cdots y_px_2y_{p+1}\cdots y_{n-2}\quad{\rm and}\quad
v_q=x_2y_1\cdots y_qx_1y_{q+1}\cdots y_{n-2},~
$$
where $0\le p,q,~p+q=n-2$.

Suppose that
$$
f=\sum_p\lambda_pw_p+\sum_q\mu_qv_q
$$
and that at least one of the coefficients $\lambda_p$ is nonzero. We may also assume that $\mu_0=0$
if $\lambda_0\ne 0$. If all $\lambda_p=0$ for  $p$ even and all $\mu_q=0$ for $q$ even, then 
$f\in Id^*(A_T)\cap Id^*(\widetilde A_T)$.

 Denote
$$
t=\max \{p|p~{\rm even~~and}~\lambda_p\ne 0\}.
$$

Then there exists odd $j$ such that $j+t=2T+1$. Hence
$$
\varphi(f)=\lambda_tz_ja^tba^m=\lambda_tz_{t+1}\ne 0\quad {\rm in}~~A_T
$$
for the evaluation $\varphi$ such that $\varphi(x_1)=z_j,\varphi(x_2)=b,
\varphi(y_1)=\cdots=\varphi(y_{n-2})=a$. 

Similarly,
$$
\widetilde\varphi(f)=\lambda_tz_{m+1}^2 \quad {\rm in}~~\widetilde A_T
$$
if $$\widetilde\varphi(x_1)=z^1_j, \widetilde\varphi(x_2)=b_1,
\widetilde\varphi(y_1)=\cdots=\widetilde\varphi(y_{n-2})=a.$$

  It follows that
$$Id^*(A_T)\cap P^*_n=Id^*(\widetilde A_T)\cap P^*_n,$$ provided that $n\le 2T$.
\hfill $\Box$

%\vskip .1in

\begin{remark}\label{r}
%{\bf Remark 1}. 
%
It follows from Lemma \ref{l1}, Lemma \ref{l2}, and Lemma \ref{l6}, that $*$-codimensions
of small degree of $\widetilde A_T$ are polynomially bounded, $$c^*_n(\widetilde A_T)\le n^3
\
\hbox{ if}
\ \ 
n\le 2T.$$ Also, any multilinear $*$-identitiy of $\widetilde A_T$ of degree $n\le 2T$ is an identity
of all $\widetilde A_{T+1},\widetilde A_{T+2},\ldots~~$.
\end{remark}
%\vskip .1in

Unlike  $A_T$, algebra $\widetilde A_T$ has an overexponential $*$-codimension growth.

\begin{lemma}\label{l9}
Let $n\ge 4T+3$. Then
\begin{equation}\label{ea1}
c^*_n(\widetilde A_T)>\left[\frac{n}{2T+1}-1\right]!,
\end{equation}
where $[t]$ denotes the integer part of real number $t>0$.
\end{lemma}

{\em Proof}.
Denote
$$
w_\sigma=x_0y_1\cdots y_{2T}x_{\sigma(1)}y_{2T+1}\cdots y_{4T}x_{\sigma(2)}\cdots
x_{\sigma(m)}y_{2mT+1}\cdots y_{2mT+j},
$$
where
$\sigma\in S_m, 0\le j\le 2T$. Since
$$
z^1_1a^{2T}b_1a^{2T}\cdots a^{2T}b_ma^j=z^{m+1}_{j+1}\ne 0,
$$
while
$$
z^1_1a^{2T}b_{\sigma(1)}a^{2T}\cdots a^{2T}b_{\sigma(m)}a^j= 0,
$$
for any $e\ne\sigma\in S_m$, all monomials $w_\sigma$ of degree $n=(2T+1)m+j+1$
are linearly independent modulo $Id^*(\widetilde A_T)$. 

Hence
\begin{equation}\label{ea2}
c^*_n(\widetilde A_T)\ge c^*_{m+1,n-m-1}(\widetilde A_T)\ge m!~~.
\end{equation}

Since $$(2T+1)m=n-j-1\ge n-(2T+1),$$ we have $$m\ge\frac{n}{2T+1}-1$$ and (\ref{ea2}) yields
 inequality (\ref{ea1}).
\hfill $\Box$
\vskip .2in

Now, let $\Phi[Z]$ be the polynomial ring over $\Phi$ and let $\Phi[Z]_0$ be its subring
of polynomials with the zero constant term. Given an integer $N\ge 1$, denote by $R_N$ the quotient
$$
R_N=\frac{\Phi[Z]_0}{(Z)^{N+1}},
$$
where $(Z)^{N+1}$ is the ideal of $\Phi[Z]_0$ generated by $Z^{N+1}$.

 Denote 
$B(T,N)=\widetilde A_T\otimes R_N$. Then
\begin{equation}\label{ea2a}
P^*_{k,n-k}(B(T,N))=P^*_{k,n-k}(\widetilde A_T),
\
\hbox{for all}
\
0\le k\le n\le N,
\end{equation}
whereas
\begin{equation}\label{ea3}
P^*_{k,n-k}(B(T,N))=0,
\
\hbox{for all}
\
n\ge N+1.
\end{equation}

Given two infinite series of integers $T_1,T_2,\ldots$ and $N_1,N_2,\ldots$ such that
$$
0<T_1<N_1<\ldots<T_j<N_j<\ldots,
$$
we define an algebra $C(T_1,T_2,\ldots,N_1,N_2,\ldots)$ as the direct sum
$$
C(T_1,T_2,\ldots,N_1,N_2,\ldots)=B(T_1,N_1)\oplus B(T_2,N_2)\oplus\cdots~~.
$$

The next statement easily follows from Lemma \ref{l2}, Lemma \ref{l8}, and relations
(\ref{ea2a}), (\ref{ea3}).

\begin{lemma}\label{l10}
Let $C=C(T_1,\cdots,N_1,\cdots)$. Then
\begin{itemize}
\item[$\bullet$]
$c^*_n(C)=c^*_n(\widetilde A_{T_1}),$ for all $n\le N_1$;
\item[$\bullet$]
$c^*_n(C)=c^*_n(\widetilde A_{T_j}),$ for all $j\ge 2, N_{j-1}+1\le n\le T_j$ ;
\item[$\bullet$]
$c^*_n(\widetilde A_{T_j})\le c^*_n(C)\le c^*_n(\widetilde A_{T_{j}})+c^*(\widetilde A_{T_{j+1}}),$
for all $j\ge 2, T_j<n\le N_j$.
\end{itemize}
\end{lemma}
\hfill $\Box$

\begin{lemma}\label{l11}
Let $C=C(T_1,\ldots,N_1,\ldots)$. Then $c^*_n(C)\le 3nc^*_{n-1}(C)$.
\end{lemma}
{\em Proof}.
Fix $n\ge 3$ and $1\le k\le n-1$. Denote by $f_1,\ldots,f_m$ a basis of $P_{k,n-k-1}^*$ modulo
$Id^*(C)$, where $f_j,1\le j\le m$, are monomials in $x_1,\ldots,x_k,y_1,\ldots,y_{n-k-1}$ and
$m=c^*_{k,n-k-1}$. 
Denote also by $g_1,\ldots,g_t$ a basis consisting of monomials in $x_1,\ldots,x_{k-1}$,
$y_1,\ldots,y_{n-k}$ of $P^*_{k-1,n-k}$ modulo $Id^*(C)$, $t=c^*_{k-1,n-k}(C)$. 

Then modulo $Id^*(C)$, 
the subspace $P^*_{k,n-k}$ coincides with the span of products
$$
f^i_1y_i,\ldots,f^i_my_i,g^j_1x_j,\ldots,g^j_tx_j,~1\le i\le n-k, 1\le j\le k,
$$
where
$$
\begin{array}{c}
f^i_p=f_p(x_1,\ldots,x_k,y_1,\ldots,y_{i-1},y_{i+1},\ldots,y_{n-k}),
\\
\\
g^j_q=g_q(x_1,\ldots,x_{j-1},x_{j+1},\ldots,x_k,y_1,\ldots,y_{n-k}).
   \end{array}
$$

Hence
\begin{equation}\label{ea4}
c^*_{k,n-k}(C)\le n(c^*_{k-1,n-k}(C)+c^*_{k,n-k-1}(C)).
\end{equation}

It follows from (\ref{ea4}) and the next  inequalities
$$
{n\choose k}\le n{n-1\choose k},~~{n\choose k}\le n{n-1\choose k-1}
$$ 
 that
\begin{equation}\label{ea5}
{n\choose k}c^*_{k,n-k}(C)\le n\left[{n-1\choose k-1}c^*_{k-1,n-k}(C)+
{n-1\choose k}c^*_{k,n-k-1}(C)\right].
\end{equation}

Inequality (\ref{ea5}) implies that
$$
\sum_{k=1}^{n-1}{n\choose k}c^{*}_{k,n-k}(C)\le 2\sum_{j=0}^{n-1}{n\choose j-1}
c^*_{j,n-j-1}(C)=2nc^*_{n-1}(C).
$$

Finally, since $c^*_{0,n}=1$ and $c^*_{n,0}=1$ for $n\ge  3$, we have
$$
c^*_n(C)\le 3n c^*_{n-1}(C).
$$\hfill $\Box$

We are now ready to construct a family of examples of algebras with involution without
$*$-PI-exponent. The following is the third main result of this paper.

\begin{theorem}\label{t3}
For any real number $\alpha>1$, there exists an algebra $C_\alpha$ such that
$$
\underline{{\rm exp}}^*(C_\alpha)=1,~~\overline{{\rm exp}}^*(C_\alpha)=\alpha.
$$
\end{theorem}
{\em Proof}.
Given $\alpha>1$, we construct an algebra $C_\alpha$ as $C(T_1,\ldots, N_1,\ldots)$
by the special choice of the sequences $T_1,T_2,\ldots$ and $N_1,N_2,\ldots$~.

First, we fix $T_1$ such that $n^3<\alpha^n$, for all $n\ge T_1$. By Lemmas \ref{l1}, \ref{l8}
and \ref{l9}, there exists $N_1$ such that
$$
\left\{
               \begin{array}{l}
        c^*_n(\widetilde A_T)<\alpha^n \ \hbox{if}\ n=N_1-1   \\
         c^*_n(\widetilde A_T)\ge\alpha^n\ \hbox{if}\ n=N_1.
               \end{array}
\right.
$$

Then by Lemma \ref{l10} and Lemma \ref{l11},
$$
\alpha^n\le c^*_n(C)\le 3n\alpha^n~~\hbox{if}~~n=N_1.
$$

On the other hand, $c^*_{N_1+1}\le (N_1+1)^3$ by the choice of $N_1$. We now set $T_2=2N_1$.

Suppose that $T_1,N_1,\ldots,T_{k-1},N_{k-1}, T_k$ have already been choosen. Then as before, applying Lemmas \ref{l1},
\ref{l8}, \ref{l9} and \ref{l10}, one can find $N_k$ such that
\begin{equation}\label{ea5a}
\left\{
               \begin{array}{l}
        c^*_n(C)<\alpha^n\ \hbox{if}\ n=N_k-1   \\
         c^*_n(C)\ge\alpha^n\ \hbox{if}\ n=N_k.
               \end{array}
\right.
\end{equation}

Moreover,
\begin{equation}\label{ea6}
\left\{
               \begin{array}{l}
        c^*_n(C)\le 3n\alpha^n   \\
         c^*_{n+1}(C)\le (n+1)^3
               \end{array}
\right.
\end{equation}
if $n=N_k$.

Denote by $C_\alpha$ the obtained algebra $C(T_1,\ldots,N_1,\ldots)$. Since $c^*_n(C_\alpha)\ne 0$
for all $n\ge 1$, relations (\ref{ea5a}), (\ref{ea6}) give us the equations
$$
\underline{{\rm exp}}^*(C_\alpha)=1,~~\overline{{\rm exp}}^*(C_\alpha)=\alpha
$$
and we have thus completed the proof.
\hfill $\Box$

\subsection*{Acknowledgements}
Repov\v s was supported by the Slovenian Research Agency program P1-0292
and grants N1-0278, N1-0114 and N1-0083.
Zaicev was supported by the  Russian Science Foundation grant 22-11-00052.

\end{document}